\documentclass[12pt,a4paper]{article}
\usepackage[OT1]{fontenc}
\usepackage[latin1]{inputenc}
\usepackage[english]{babel}
\usepackage{amsfonts}
\usepackage{amsthm}
\usepackage{amsmath}
\newtheorem{teo}{Theorem}[section]
\newtheorem{prop}[teo]{Proposition}
\newtheorem{lem}[teo]{Lemma}
\newtheorem{coro}[teo]{Corollary}
\newtheorem{defi}[teo]{Definition}
\theoremstyle{definition}
\newtheorem{rem}[teo]{Remark}
\newtheorem{ejem}[teo]{Example}

\def\M{{\cal M}}
\def\um{{\cal U}_{\cal M}}

\def\N{{\cal N}}
\def\la{{\lambda}}

\def\lpm{{L}^p({\cal M})}
\def\lpn{{L}^p({\cal N})}

\def\beah{{\cal B}_2({\cal H})_{sh}}

\def\u2{U_2({\cal H})}

\def\o{ {\cal O} }
\def\f{ {\cal F} }
\def\g{ {\cal G} }
\def\ad{ {\mbox{ad\,}} }
\def\noi{ \noindent }

\begin{document}

\title{\vspace*{0cm}Homogeneous manifolds from non commutative measure spaces}
\date{}
\author{Esteban Andruchow, Eduardo Chiumiento and Gabriel Larotonda\footnote{All authors partially supported by Instituto Argentino de Matemática and CONICET.}}

\maketitle

\abstract{\footnotesize{\noindent
Let $\M$ be a finite von Neumann algebra with a faithful trace $\tau$. In this paper we study metric geometry of homogeneous spaces $\o$ of the unitary group $\um$ of $\M$, endowed with a  Finsler quotient metric induced by the $p$-norms of $\tau$, $\|x\|_p=\tau(|x|^p)^{1/p}$, $p\geq 1$. The main results include the following. The unitary group carries on a rectifiable  distance $d_p$ induced by measuring the length of curves with the $p$-norm. If we identify $\o$ as a quotient of groups, then there is a natural quotient distance $\dot{d}_p$ that metrizes the quotient topology. On the other hand, the Finsler quotient metric defined in $\o$ provides a way to measure curves, and therefore, there is an associated rectifiable distance $d_{\o , p}$\,. We prove that the distances $\dot{d}_p$ and $d_{\o , p}$ coincide. Based on this fact, we show that the metric space $(\o,\dot{d}_p)$ is a complete path metric space. The other problem treated  in this article is the existence of metric geodesics, or  curves of minimal length,  in $\o$. We give two abstract partial results in this direction. The first concerns  the initial values problem and the second  the fixed endpoints problem. We show how these results apply to several examples. In the process, we improve some results about the metric geometry of $\um$ with the $p$-norm.}}


\section{Introduction}

In this paper we study metric properties of smooth homogeneous spaces $\o$ of the unitary group ${\cal U}_{\cal M}$ of a finite von Neumann algebra $\M$. If $\tau$ is a finite faithful trace on $\M$, the $p$-norms on $\M$ ($p\ge 1)$ induced by the trace,
can be used to turn $\um$ into a complete metric space. This is achieved by giving $T\um$ the Finsler metric that is given by the $p$-norm at any point of $\um$, recalling that its Lie algebra can be identified with $\M_{sh}$, the skew-hermitian elements of $\M$.
The $p$-norms in $\M$ can be used to turn $\o$ into a metric space, in two different ways. First, measuring the distance $\dot{d}_p$ between classes. Second, with the rectifiable distance $d_p$ induced by the Finsler metric just mentioned.


In both cases it is necessary to assume that the isotropy groups $G_x\subset \um$ are closed in the $p$-norm (or else one ends with a pseudo-distance where $d(x,y)=0$ might not imply $x=y$ in $\o$). Note that when $G_x$ is the unitary group of a von Neumann subalgebra of $\M$, then $G_x$ is $p$-norm closed in $\um$. The first metric is well-known (as it is the quotient metric on a quotient of a topological metricable group), and moreover it is known that $(\o,\dot{d}_p)$ is a complete metric space \cite[p.109]{takesaki3}. We show here that for $p>1$, these two metrics coincide with no additional hypothesis (Theorem \ref{iguales}). 

Next we study the existence of metric geodesics, or  short paths, for the given metric.  We give a first step in that direction in Theorem \ref{geod}, where it is shown that, under suitable hypothesis (in particular $p$ even), the curves of the form $\delta(t)=e^{tz}\cdot x$, with \textit{minimal symbol} $z\in \M_{sh}$, that is $\|z\|_p\le \|z-y\|_p \;\mbox{ for any }y\in \g_x$ are minimizing in $\o$ up to a critical $t$, among a certain family of rectifiable curves.

We also  show (Theorem \ref{geodesico}) that there exists a certain set in $\o$ containing $x$ such that any point there can be joined to $x$ with such a curve, which is shorter than any other smooth curve joining the same endpoints which does not leave the mentioned set.

In some  examples, for instance quotient spaces of $\um$ with the unitary group ${\cal U}_{\cal N}$, where ${\cal N}$ is a von Neumann subalgebra of the center of $\M$, this set is in fact an open uniform neighborhood of $x$ in $\o={\cal U}_{\cal M}/{\cal U}_{\cal N}$.

This paper is organized as follows: in Section \ref{section2}, we introduce the necessary definitions (finite von Neumann algebras,  smooth homogeneous structures, rectifiable distances and uniform convexity in Banach spaces). Section \ref{section3} contains known results on the rectifiable distance in the unitary group $\um$ of a finite von Neumann algebra $\M$, in particular minimality of geodesics and local convexity of the geodesic distance. In Section \ref{section4} we introduce the rectifiable distance in the homogeneous space $\o$, and we prove several facts on the metric and topological properties of this space, including the coincidence between the quotient topology and the topology induced by the rectifiable distance in $\o$. These imply the completeness of the metric space $\o$ with the rectifiable distance. In Section \ref{section5} we treat the problem of minimality of geodesics in $\o$, and  prove two partial results related to the initial values problem and the fixed endpoints problem. We finish the paper with a collection of examples where our results apply, among them the homogeneous spaces $\o={\cal U}_{\cal M}/{\cal U}_{\cal N}$, where ${\cal N}$ is a von Neumann subalgebra of the center of $\M$.

\section{Definitions and background}\label{section2}

Let $\M$ be a von Neumann algebra with a finite and faithful trace $\tau$. Denote by $\M^{\times}$ and $\um$ the
groups of invertible and unitary operators of $\M$. Let $1\le p\le \infty$, and denote with $\lpm$ the noncommutative $L^p$ space of $(\M,\tau)$, that is, the completion of $\M$ relative to the norm $\|\cdot\|_p$, where $\|x\|_p^p=\tau (|x|^p)$. When $p=2$, $\lpm$ is a Hilbert space with the inner product $<a,b>_{\tau}=\tau(b^*a)$. We use the subscript $h$ (resp. $sh$) to denote the sets of hermitian (resp. skew-hermitian) operators. The symbol $\| \ \|$ denotes the usual operator norm of $\M$.

If $f:X\to Y$ is a smooth map between manifolds, we will use $f_*:TX\to TY$ to denote the differential of $f$ and $f_{*x}:T_xX\to T_{f(x)}Y$ its specialization. Let $\o$ be a topological space on which $\um$ acts continuously and transitively, such that for any element $x\in\o$, the subgroup $G_x=\{u\in\um: u\cdot x=x\}$, called the {\it isotropy group} of the action at $x$, is a closed submanifold of $\um$. This implies that $\o$ can be endowed with a differentiable manifold structure, in a way such that the map
$$
\pi=\pi_{x}: \um \to \o , \ \ \pi_{x}(u)=u\cdot x
$$
is a smooth submersion. Therefore $\o$ is  a homogeneous space of the group $\um$.

\begin{rem}\label{parte}
For $x\in \o$, denote by $\g_x$ the Banach-Lie algebra of the isotropy subgroup $G_x$. Since we are assuming that $\pi$ is a smooth submersion, and that $\o$ is given the differentiable structure that induces the final topology on $T_x\o$, there exists a closed supplement ${\cal F}_x\subset \M_{sh}$ such that $\M_{sh}\simeq \g_x\oplus {\cal F}_x$, and a smooth section $s_x:T_x\o\to {\cal F}_x$, $(\pi_x)_{*1}\circ s_x=id_{T_x\o}$, $s_x\circ (\pi_x)_{*1}=P_{{\cal F}_x}$, where the last expression denotes the unique bounded projection in ${\cal B}(\M_{sh})$ with rank ${\cal F}_x$ and kernel $\g_x$. The tangent space $T_x\o$ can be normed with the (uniform)norm,
$$
\|V\|_x=\inf\{\|z-y\|: y\in \g_x\},
$$
where $z\in \M_{sh}$ is any lift of $V$, i.e. $\pi_*(z)=V$. Note that
$$
\|s_x(V)\|=\|s_x((\pi_x)_{*1}(z))\|=\|P_{{\cal F}_x}(z)\|=\|P_{{\cal F}_x}(z-y)\|\le\|P_{{\cal F}_x}\|\|z-y\|
$$
for any lift $z\in \M_{sh}$ of $V\in T_x\o$ and any $s\in\g_x$. Thus $\|s_x(V)\|\le \|P_{{\cal F}_x}\|\|V\|_x$. The norm of $P_{{\cal F}_x}$ does not depend on the point $x \in \o$ since
$$ \|P_{{\cal F}_{u \, \cdot  \,  x}}\|=\| Ad_u \circ P_{{\cal F}_{x}} \circ Ad_{u^*}\|=\|P_{{\cal F}_x}\|, $$
where $Ad_u:\M \longrightarrow \M$, $Ad_u(z)=uzu^*$.
Therefore, there exists a constant $C_{\o}$ depending only on the differentiable structure such that
$$
\|s_y(V)\|\le C_{\o} \|V\|_y\;\mbox{ for any }y\in \o \;\mbox{ and any }V\in T_y\o.
$$
\end{rem}

\subsection{Quotient metrics}

Throughout this article $1< p< \infty$ unless otherwise stated, and $L_p$ denotes the length functional for piecewise smooth curves in $\um$, measured with the $p$-norm:
$$
L_p(\alpha)=\int_{t_0}^{t_1} \|\dot{\alpha}(t)\|_p\, dt,
$$
while (unless otherwise stated) smooth means $C^1$ and with nonzero derivative, relative to the uniform topology of $\M$. Let us introduce some notation. The action of $\um$ on $\o$ induces two kind of maps. If one fixes $x\in\o$, one has the submersion
$$
\pi_x:\um\to \o, \ \pi_x(u)=u\cdot x,\; u\in\um.
$$
If one fixes $u\in\um$ one has the diffeomorphism
$$
\ell_u:\o\to \o ,  \ \ell_u(x)=u\cdot x, \; x\in\o.
$$
If $x\in\o$ and $X\in T_x\o$,  put
$$
\|X\|_{x,p}=\inf\{ \|z\|_p: z\in \M_{sh}, (\pi_x)_{*1}(z)=X\}.
$$
We call this Finsler metric the \textit{quotient metric} of $\o$, because it is the quotient metric in the metric linear space $T_x\o$ if one identifies it with $\M_{sh}/\g_x$.  Indeed, since $\g_x=\ker (\pi_x)_{*1}$, if $z\in \M_{sh}$ with $(\pi_x)_{*1}(z)=X$, then
$$
\|X\|_{x,p}=\inf \{ \|z-y\|_p: y\in\g_x\}.
$$
We shall omit the index $p$ since it is fixed in any discussion. One of the main features of this metric in $\o$ is that it is invariant by the group action (or in other words, that the group acts isometrically on the tangent spaces): a straightforward computation shows that if $x\in\o$, $X\in T_x\o$ and $u\in\um$, $\|(\ell_u)_{*x}(X)\|_{u\cdot x}=\|X\|_x$. Note that when $p=\infty$, this is the metric that arises naturally if we regard $\o$ as an homogeneous space of $\um$, as a Banach-Lie group for the topology induced in $\M_{sh}$ by the uniform norm, as discussed in Remark \ref{parte}.

\begin{rem}
Since the action is transitive, we shall drop the index $x\in\o$ for the maps involved, when there is no possibility of confusion. The isotropy group will be denoted by $G$ and the Lie algebra by $\g$, and $\pi:\um\to \o$ will denote the smooth  submersion.
\end{rem}

\begin{rem}\label{clark}
Recall Clarkson's inequalities \cite{kosaki} for $\lpm$ spaces. Let $a,b\in \lpm$,  $1/p+1/q=1$, then
$$
(\|a+b\|_p^q+\|a-b\|_p^q)^{\frac1q}\le 2^{\frac1q} (\|a\|_p^p+\|b\|_p^p)^{\frac1p} \mbox{ if } 1< p\le 2, \mbox{ and }
$$
$$
(\|a+b\|_p^p+\|a-b\|_p^p)^{\frac1p}\le 2^{\frac1q} (\|a\|_p^p+\|b\|_p^p)^{\frac1p} \mbox{ if } 2\le p<\infty.
$$
From them it can be easily derived that $\lpm$ is uniformly convex and uniformly smooth, and that for any convex closed set $S\subset \lpm_{sh}$ there exists a continuous map $Q_{S,p}:\lpm_{sh}\to S$ which sends $x\in \lpm_{sh}$ to its best approximant $Q_{S,p}(x)\in S$, i.e.
$$
\|x-Q_{S,p}(x)\|_p\le \|x-s\|_p
$$
for any $s\in S$. The map $Q_{S,p}$ is single-valued and continuous since $\lpm$ is uniformly convex and uniformly smooth (see for instance \cite{chongli}). Omitting the index $p$ for convenience, note that
$$
\|Q_S(x)\|_p\le \|Q_S(x)-x\|_p+\|x\|_p\le \|0-x\|_p+\|x\|_p=2\|x\|_p
$$
and also that
$$
\|x-Q_S(x)-s\|_p\ge \|x-Q_S(x)\|_p
$$
for any $s\in S$, hence $Q_S(x-Q_S(x))=0$, namely $Q_S\circ (1-Q_S)=0$. Also, for any $\lambda\in\mathbb{R}$, $Q_S(\lambda x)=\lambda Q_S(x)$. Calling $\bar{Q}_S=1-Q_S$, we have
$$
S=\bar{Q}_S ^{-1}(0)=Im(Q_S),\qquad Q_S^{-1}(0)=Im(\bar{Q}_S),\qquad
$$
and also
$$
\bar{Q}_S ^2=\bar{Q},\quad Q_S ^2=Q_S,\quad \bar{Q}_S\circ Q_S=Q_S\circ \bar{Q}_S=0,
$$
which shows that $Q_S$ has some of the  properties of the linear projection (when $p=2$).
\end{rem}

Let $x\in {\cal O}$,  $G$ be the isotropy group, ${\cal G}$ the Lie algebra of $G$ and $\overline{\cal G}^p$ its closure in $\lpm_{sh}$. Let $Q=Q_{\g}$ be the projection to the best approximant in $\overline{\g}^p$. Let
$$
\g^{\perp_p}=Q^{-1}(0)=\{x\in \lpm_{sh}:\|x\|_p\le \|x-y\|_p \; \mbox{ for any } \; y\in \g\}.
$$
Then any element $z\in \lpm_{sh}$ can be uniquely decomposed as
$$
z=z-Q(z)+Q(z),
$$
where $z-Q(z)=(1-Q)(z)\in \g^{\perp_p}$ and $Q(z)\in \overline{\g}^p$.

\begin{rem}
In particular, for $1< p<\infty$, the quotient metric of $\o$ is given by
$$
\|X\|_x:=\|z-Q(z)\|_p,
$$
where $z\in \M_{sh}$ is any element such that $\pi_{*1}(z)=X$. Note that there always exists such $z$ since $\pi_{*1}$ is surjective. We call $z_0=z-Q(z)\in\lpm_{sh}$ a \textit{minimal lifting} of $X$. A word of caution: the map $Q$ depends on the chosen parameter $p>1$.

\bigskip

Note that if $p=2$, this metric is Riemannian. Indeed, if $Q_x=1-P_x$ is the orthogonal projection onto $\overline{\g_{x}}^p$, then each $z\in \M_{sh}$ can be uniquely decomposed as $$z=z-Q_x(z)+Q_x(z)=z_0+Q_x(z)$$ and $z_0=P_x(z)$ is orthogonal to $\g_{x}$ hence
$$
\|z-y\|_2^2=\|z_0+Q_x(z)-y\|_2^2=\|z_0\|_2^2+\|Q_x(z)-y\|_2^2\ge \|z_0\|_2^2
$$
for any $y\in \g_{2,x}$, which shows that
$$
\|X\|_x=\inf \{ \|z-y\|_2: y\in\g_{x}\}=\|z_0\|_2
$$
where $z_0$ is the unique vector in $\g_{x}^{\perp}$ such that $(\pi_x)_{*1}(z_0)=X$.

\bigskip

We shall denote with $\overline{T_x\o}^p$ the completion of $T_x\o$ relative to the $p$-quotient metric. Then $\pi_{*1}$ extends naturally to a surjective linear map $\pi^p_*:\lpm_{sh}\to \overline{T_x\o}^p$, since
$$
\|\pi_*(y_n)-\pi_*(z_n)\|_x=\|y_n-z_n-Q(y_n-z_n)\|_p\le \|y_n-z_n\|_p
$$
and then one can put $\pi^p_*(z_0)=\lim_n\pi_{*1}(z_n)$ disregarding the particular sequence $(z_n)_{n  \geq 1}$ such that $z_n\to z_0\in \lpm_{sh}$.
\end{rem}

\begin{lem}\label{perpep}
Let $p$ be an even positive integer. Let $x\in\o$ and $X\in \overline{T_x\o}^p$.  An element $z_0\in\lpm_{sh}$ with $\pi_*^p(z_0)=X$ is a minimal lifting for $X$ if and only if $\tau(z_0^{p-1}y)=0$ for all $y\in\g$. For any $X\in \overline{T_x\o}^p$ there exists a unique minimal lifting $z_0\in \g^{\perp_p}$ such that $\|z_0\|_p=\|X\|_x$.
\end{lem}
\begin{proof}
The proof is straightforward, see for instance  \cite[Lemma 4.3]{alr}.
\end{proof}

\section{Metric structure of $\um$}\label{section3}

In this section we recall and complete certain facts from \cite{duranmatarecht}, concerning the minimality of geodesics in $\um$, and in addition we prove a local convexity result. The following elementary lemma will be used in the proof of Theorem \ref{minimalidadunitarios}, its proof can be found in \cite[Lemma 3.4]{alr}.

\begin{lem}\label{fseg}
Let $C,\varepsilon >0$, let $f(-\varepsilon,1+\varepsilon)\to \mathbb R$ be a non constant real analytic function such that $f'(s)^2\le C f''(s)$ for any $s\in [0,1]$. Then $f$ is strictly convex in $(0,1)$.
\end{lem}

\begin{rem}\label{remarko}
\begin{enumerate}
\item The map $\exp(x)=e^x$, $exp: \M_{sh}\to \um$ is surjective.
\item The exponential map is a diffeomorphism between the sets
$$
\M_{sh}\supset \{z\in\M_{sh}: \|z\|<\pi\}\to \{u\in\um: \|1-u\|<2\}.
$$
\item
Moreover, $exp:\{z\in \M_{sh}: \|z\|\le \pi\}\to \um,$ is surjective.
\end{enumerate}
\end{rem}

\medskip

For $a,b\in \M$, let $R_a,L_a:\M\to \M$ stand respectively for the right and left multiplication,  and let $\ad a=R_a-L_a:\M\to \M$ stand for the adjoint operator. Then in \cite[Lemma 3.3]{alr} it is proved the following lemma. Its proof can be adapted to our context without any modification.
\begin{lem}\label{expo}
Let $a,b\in \M$. Then
$$
exp_{*a}(b)=\int\limits_0^1 e^{(1-t)a}be^{ta}\, dt=e^a\,F(\ad a)b= F(\ad a)(e^a\, b),
$$
where $F(z)=\frac{e^z-1}{z}=\sum_{n\ge 0}\frac{z^{n}}{(n+1)!}$. The differential is invertible at $a$ if and only if $\sigma(ad(a))\cap \{2k\pi i\}=\emptyset$ ($k\ne 0$), and then $
exp_{*a}^{-1}(w)=e^{-a}F(\ad a)^{-1}w$. In particular if $\|a\|<\pi$ then $exp_{*a}$ is invertible. If  $a\in \M_{sh}$, the differential is a contraction, that is $\|exp_{*a}(b)\|_p\le \|b\|_p$ for any $p\in [1,\infty]$.
\end{lem}

\begin{rem}\label{bilinear}
Let $a,b,c\in \M_{sh}$, let $H_a:\M_{sh}\to \mathbb R$ stand for the symmetric bilinear form given by
$$
H_a(b,c)=(-1)^{\frac{p}{2}}p \sum_{k=0}^{p-2}\, \tau(a^{p-2-k}b a^k c).
$$
If $Q$ is the quadratic form associated to $H$, then (\cite[Lemma 4.1]{convexg} and equation (3.1) in \cite{cocomata}):
\begin{enumerate}
\item $Q_a([b,a])\le 4 \|a\|^2 Q_a(b)$.
\item $Q_a(b)=p\|ba^{\frac{p}{2}-1}\|_2^2+\frac{p}{2}\sum_{l+m=n-2}\|a^l (ab+ba)a^m\|_2^2$.
\end{enumerate}
In particular $H_a$ is positive definite for any $a\in \M_{sh}$.
\end{rem}

\noi The following theorem collects several results concerning the rectifiable $p$-distance in the unitary group of $\M$, such as: minimality of geodesics, uniqueness of such geodesics, comparison with the usual $p$-distance, and finally a fundamental convexity result which improves the one stated in \cite{convexg}.

\begin{teo}\label{minimalidadunitarios}
Let $2\le p<\infty$. The following facts hold.
\begin{enumerate}
\item
Let $u\in\um$ and $x\in\M_{sh}$ with $\|x\|\le \pi$. Then the curve
$\mu(t)=ue^{tx}$, $t\in[0,1]$ is shorter than any other smooth curve in
$\um$ joining the same endpoints, when we measure them with the length functional $L_p$. Moreover, if $\|x\|<\pi$, the curve $\mu$ is unique with this property among all the $C^2$ curves in $\um$.
\item
Let $u_0, u_1 \in\um$. Then there exists a minimal geodesic curve joining them.
If $\|u_0-u_1\|<2$,  this geodesic is unique among all the $C^2$ curves there.
\item
The diameter of $\um$ is $\pi$ for all the $p$-norms.
\item
If $u,v\in\um$ then
$$
\sqrt{1-\frac{\pi^2}{12}} \; d_p(u,v) \le \|u-v\|_p\le d_p(u,v).
$$
In particular the metric space $(\um, d_p)$ is complete.
\item Let $p$ be an even positive number, $u,v,w,\in \um$, with
$$
\|u-v\|<\sqrt{2},\qquad \|w-v\|<\sqrt{2}-\|u-v\|.
$$
Let $\beta$ be a short geodesic joining $v$ to $w$ in $\um$. Then the rectifiable $p$-distance between $u$ and $\beta$ is a strictly convex function, provided $u$ does not belong to any prolongation of $\beta$.
\end{enumerate}
\end{teo}
\begin{proof}
The minimality was proved in \cite[Theorem 5.4]{otrococo}.
Let us prove that if $\|x\|<\pi$, then $\mu$ is unique with the minimality property among all the smooth curves.
To do this, we shall follow a standard procedure, using the first variation formula, in this case, for the functional $F_p$ which is given by
$$
F_p(\gamma)=\int_0^1 \|\dot{\gamma}(t)\|_p^p d t,
$$
if $\gamma(t)\in\um$, $t\in[0,1]$. Let $\gamma_s(t)$, $t\in [0,1]$, $s\in (-r,r)$ be a $C^2$ variation of the curve $\gamma$, i.e.
$\gamma_s(t)\in \um$, for all $s,t$, the map $(s,t)\mapsto \gamma_s(t)$ is $C^2$ and $\gamma_0(t)=\gamma(t)$. We shall use a formula for $\frac{d}{d s} F_p(\gamma_s)\arrowvert_{s=0}$, obtained in \cite{convexg} in the context of a $C^*$-algebra with trace. As in classical differential geometry, we shall call the expression obtained the first variation formula.  Let
$$
V_s=\frac{d}{d t}{\gamma}_s \hbox{ and } W_s=\frac{d}{d s}\gamma_s.
$$
With lower case types we denote the left translations $v_s=\gamma_s^* V_s$ and $w_s=\gamma_s^*W_s$. Note that $V_s, W_s\in (T\um)_{\gamma_s}$ whereas $v_s, w_s\in \M_{sh}$.  Then
$$
\frac{(-1)^{p/2}}{p} \frac{d}{d s} F_p(\gamma_s)=\tau(v_s^{p-1}w_s)\arrowvert_{t=0}^{t=1}-\int_0^1 \tau( \frac{d}{d t}[v_s^{p-1}] w_s ) d t.
$$
Suppose that $\gamma(t)\in\um$ is a $C^2$ minimal curve, and let $\gamma_s(t)$ be a variation, with fixed endpoints $\gamma(0)$ and $\gamma(1)$, i.e. $\gamma_s(0)=\gamma(0)$ and $\gamma_s(1)=\gamma(1)$ for all $s$. Then
$\frac{d}{ds}F_p(\gamma_s)|_{s=0}=0$, and thus
$$
0=\tau(v_0^{p-1}w_0)|_{t=0}^{t=1}- \int_0^1 \tau(w_0 \frac{d}{dt}(v_0^{p-1})) d t.
$$
The fixed endpoints hypothesis implies that the first term vanishes. Then
$$
\int_0^1 \tau(w_0 \frac{d}{dt}(v_0^{p-1})) d=0
$$
for any variation $\gamma_s$ with fixed endpoints.
Let us denote by $Z(t)=\frac{d}{dt}(v_0^{p-1})$ and by $A(t)=w_0(t)$. Both $A$ and $Z$ are continuous fields in $\M_{sh}$. The variation formula implies that
$$
\int_0^1 \tau(A(t)Z(t)) d t = 0
$$
for any continuous field  $A$ in $\M_{sh}$ such that $A(0)=A(1)=0$.  We claim that this condition implies that $Z(t)=0$ for all $t$.

First note that the requirement that the field $A$ vanishes at $0$ and $1$ can be removed: let $f_r(t)$ be a real function which is constant and equal to $1$ in the interval $[r, 1-r]$ and such that $f(0)=f(1)=0$, with $0\le f_r(t)\le 1$ for all $t$. Let $B(t)$ be any continuous field in $\M_{sh}$ and consider  $A_r(t)=f_r(t)B(t)$. Then $\int_0^1 A_r(t) Z(t) d t =0$, and if $r\to 0$, $\int_0^1 B(t) Z(t) d t =0$. Also it is clear that the integral will vanish if $A$ is non skew-hermitian. Indeed, it is clear if $A$ is hermitian, and for general $A$,  decompose $A$ as the sum of its hermitian and skew-hermitian parts.

Consider $A(t)=-Z(t)$, then $\int_0^1 \| Z (t) \|_2^2 \, d t=0$, which implies $Z(t)\equiv 0$.
Therefore $v_0^{p-1}$ is constant, and since $p$ is even and $v_0$ is skew-hermitian, $v_0(t)=\gamma(t)^*\frac{d}{d t}\gamma(t)$ is constant, i.e. $\gamma(t)=e^{tx}$ for some $x\in \M_{sh}$.

Fact 2. It is straightforward from the first item and Remark \ref{remarko}.

Fact 3. Any pair of unitaries $u_0, u_1$ can be joined by a minimal curve of length less or equal than $\pi$. Indeed, let $x \in \M_{sh}$, $\|x \| \leq \pi$ and $e^x=u_0 ^* u_1$. Then $\mu(t)= u_0 e^{tx}$ have minimal length equals to $\| x \|_p \leq \| x \| \leq \pi$. Then the diameter is exactly $\pi$ since the unitaries $1$ and $-1$ are joined by the minimal curve $\mu(t)=e^{it\pi1}$, which has length $\pi$.

Fact 4.  Both metrics are invariant by left translation with elements of $\um$. Therefore it suffices to compare $d_p(u,1)$ and $\| u -1 \|_p$, for $u \in \um$.  Let $x=x^* \in \M$ with $\|x\| \leq \pi$ and $u=e^{ix}$. Then by item 1, $d_p(u,1)= \| x\|_p$. We follow Petz \cite{p} for the definition and properties of the spectral scale $\la_t(x)$ of $x$. It is defined by
\[ \la_t(x)= \inf \{ s \in \mathbb{R} \, : \, \tau(e_{(s\, , \, \infty)}(x) ) \leq t\, \},  \]
where $t \in (0,1)$ and $e_I(x)$ denotes the spectral projection of $x$ corresponding to an interval $I$ in $\mathbb{R}$. If $f$ is a real Borel function on $\mathbb{R}$, then by Proposition 1 in \cite{p} we have $\tau(f(x))= \int_0 ^1 f(\lambda_t(x)) \, dt.$ On the other hand, for $|s| \leq \pi$ it is easily seen that
\[   |s| (1- \frac{\pi^2}{12})^{1/2} \leq | e^{is} -1 | \leq |s|.   \]
Since $\|x\|\leq \pi$, we have $|\lambda_t(x)| \leq \pi$, $t \in [0,1]$. Then we obtain the inequality
$$
\| u -1 \|_p ^p  = \tau (|e^{ix} - 1|^p) = \int_0 ^1 |e^{i\lambda_t(x)} - 1|^p dt
 \leq \int_0 ^1 |\lambda_t(x)|^p dt = \|x \|_p ^p\,.
$$
The other inequality follows in the same fashion,
\[ \| u -1 \|_p ^p= \int_0 ^1 |e^{i\lambda_t(x)} - 1|^p dt \geq (1- \frac{\pi^2}{12})^{p/2} \int_0 ^1 |\lambda_t(x)|^p dt = \|x \|_p ^p =(1- \frac{\pi^2}{12})^{p/2} \| x \|_p ^p \,,  \]
and our claim holds.

It is  straightforward that $(\um ,d_p)$ is a complete metric space: let $(u_n)_{n  \geq 1}$ be a Cauchy sequence for $d_p$, then it is Cauchy in $\lpm$. Hence, it converges to an element $u_0 \in \lpm$. Since $(u_n)_{n  \geq 1}$ is uniformly bounded in the operator norm, it  follows that $u_0 \in \M$, and then clearly $u_0 \in \um$.

Fact 5. It was proved in \cite[Theorem 4.1]{convexg} a similar result. With a slight modification of the proof, one obtains the better convexity radius estimate stated here. We include the proof for the convenience of the reader. We may assume that $u=1$ since the action of unitary elements is isometric. Note that $\|v^*w-1\|=\|v-w\|<\sqrt{2}<2$, so we can compute $z=\log(v^*w)\in \M_{sh}$, where $\log$ indicates the principal branch of the logarithm. Let $\beta(s)=ve^{sz}$, which is a short geodesic joining $v$ to $w$ in $\um$ since $\|z\|<\pi$. Then
$$
\|1-ve^{sz}\|\le \|1-v\|+\|1-e^{sz}\|\le \|1-v\|+\|1-e^z\|=\|1-v\|+\|v-w\|<\sqrt{2},
$$
which implies that $\beta$ has an analytic logarithm,  $w_s=\log(\beta(s))=\log(ve^{sz})$, with $\|w_s\|<\frac{\pi}{2}$. Let $\gamma_s(t)=e^{tw_s}$, then $\gamma_s$ is a short geodesic joining $1$ and $\beta(s)$, of length $\|w_s\|_p=d_p(1,\beta(s))$. Then $f_p(s)=\|w_s\|_p^p=\tau((-w_s^2)^\frac{p}{2})=(-1)^{\frac{p}{2}}\tau(w_s^p)$, hence $$
f'_p(s)=
(-1)^{\frac{p}{2}}p\, Tr(w_s^{p-1} \dot{w_s})=\frac{1}{p-1} H_{w_s}(\dot{w_s},w_s),
$$
where $H$ is the bilinear form introduced in Remark \ref{bilinear}. Since $e^{w_s}=ve^{sz}$, then $e^{-w_s}\; exp_{*w_s} (\dot{w_s}) =z$ by Lemma \ref{expo}, namely
\begin{equation}\label{difexp}
z= \int_0^1 e^{-tw_s} \dot{w_s} e^{tw_s}\;dt.
\end{equation}
Thus $\tau(w_s^{p-1} \dot{w_s} )=\int_0^1 \tau(w_s^{p-1}e^{-tw_s} \dot{w_s} e^{tw_s})\;dt=\tau(zw_s^{p-1})$.
Hence
$$
f''_p(s)=
(-1)^{\frac{p}{2}}p \sum_{k=0}^{p-2}\, \tau(w_s^{p-2-k}\dot{w_s}w_s^k z)=H_{w_s}(\dot{w_s},z),
$$
and again by equation (\ref{difexp}) above, if we put $\delta_s(t)=e^{-tw_s}\dot{w_s} e^{tw_s}$, then
$$
f_p''(s)=\int_0^1 H_{w_s}(\delta_s(0),\delta_s(t))\,dt.
$$
Suppose that for this value of $s\in[0,1]$, $R_s^2:=Q_{w_s}(\dot{w_s})\ne 0$, where $Q_{w_s}$ is the quadratic form associated to $H_{w_s}$. If $K_s\subset \M_{sh}$ is the null space of $H_{w_s}$, consider the quotient space $\M_{sh}/K_s$ equipped with the inner product $H_{w_s}(\cdot,\cdot)$. An elementary computation shows that $\delta_s(t)$ lives in a sphere of radius $R_s$ of this pre-Hilbert space, hence $H_w(\delta_s(0),\delta_s(t))=R_s^2 \cos(\alpha_s(t))$,
where $\alpha_s(t)$ is the angle subtended by $\delta_s(0)$ and $\delta_s(t)$. Then, reasoning in the sphere
\begin{eqnarray}
R_s\alpha_s(t) &\le & L_0^t(\delta_s)=\int_0^t Q_{w_s}^{\frac12}(e^{-tw_s}[w_s,\dot{w_s}]e^{t w_s})\,dt\nonumber\\ &=&\int_0^t Q_{w_s}^{\frac12}([w_s,\dot{w_s}])\,dt=t\,Q_{w_s}^{\frac12}([w_s,\dot{w_s}]).\nonumber
\end{eqnarray}
By property $1.$ of Remark \ref{bilinear}, $R_s\alpha_s(t)\le t\, 2\|w_s\| R_s <R_s \pi$. So $
\cos(\alpha_s(t))\ge \cos(2t\|w_s\|)$ and then integrating with respect to the $t$-variable,
$$
f''_p(s)\ge R_s^2 \frac{\sin(2\|w_s\|)}{2\|w_s\|}>0
$$
provided $R_s\ne 0$. On the other hand, the Cauchy-Schwarz inequality for $H_{w_s}$ shows that if $R_s=0$, then
$$
(p-1)f_p'(s)=H_{w_s}(w_s,\dot{w_s})\le Q^{\frac12}_{w_s}(\dot{w_s})Q^{\frac12}_{w_s}(w_s)=0.
$$
Assume that $R_s$ is identically zero, $s\in [0,1]$. Then $f_p$ is constant with $f_p(s)=f_p(0)=\|v\|_p$ for any $s\in [0,1]$. Moreover, by property $2.$ of the remark above, $R_s=0$ implies $w_s^{\frac{p}{2}-1}z=0$ and an elementary computation involving the functional calculus of skew-hermitian operators shows that $w_s z=0$. Put $y=\log(v)$; in particular we have $yz=0$ which implies, since $w_s=\log(e^ye^{sz})$, that  $w_s=y+sz$ by the Baker-Campbell-Hausdorff formula. But since the $p$-norm of $\M_{sh}$ is strictly convex, $w_s=y+sz$ cannot have constant norm unless $y$ is a multiple of $z$, and in that case, $u$ and $\beta$ are aligned contradicting the assumption of the theorem. So there is at least one point $s_0\in  [0,1]$ where $R_{s_0}\ne 0$, so $f_p$ is non constant. There exists a positive constant $C$ such that $\frac{\sin(2\|w_s\|)}{2\|w_s\|}\ge C$, hence $f''_p(s)\ge  R_s^2 C=Q_w(\dot{w}) C$.
On the other hand $Q_w(w)=p(p-1)\|w_s\|_p^p\le p(p-1)\|w_s\|<p(p-1)\frac{\pi}{2}$, so
$$
(p-1)^2f_p'(s)^2=H_{w_s}^2(w_s,\dot{w_s})\le Q_{w_s}(\dot{w_s})Q_{w_s}(w_s)\le \frac{p(p-1)\pi}{2C} f''_p(s).
$$
and by Lemma \ref{fseg}, $f_p$ is strictly convex.
\end{proof}

\noindent The following remark justifies in part the election of the uniform topology to differentiate curves in $\um$. We show that in order to produce minimal curves, the velocity vectors of the exponential should have uniform length less or equal than $\pi$.

\begin{rem}
Let $z \in \lpm_h$ such that $\pi<\| z \| < \infty$ . By the Stone theorem the curve $\delta(t)=e^{itz}$ is $C^1$ if we differentiate it in the strong topology of the standard representation of $\M$. We claim that $\delta$ is not of minimal length  joining its endpoints in $\um$ when we measure with the $p$-norm.

In order to prove  this consider the function $f:\mathbb{R} \longrightarrow [-\pi,\pi]$ given by
\[ f(t)=\left\{
\begin{array}{ccc}
 t +2(k+1)\pi & \, \, \, \, \, \,  \, \, \, \, \, -2(k+3)\pi \leq t < -(2k +1)\pi \\
 t & \, \, \, \, \, \,  \, \, \, \, \, -\pi \leq t \leq \pi \\
 t -2(k+1)\pi & \, \, \, \, \, \, \, \, \, \, \, (2k+1)\pi < t \leq (2k+3)\pi.
\end{array}
\right.
\]
Clearly it is a Borel measurable function. Then we use the Borel functional calculus of   $z=\int_{\sigma(z)}\lambda \, de(\lambda)$ to obtain
\[ <e^{if(z)}\xi , \eta > = \int_{\sigma(z)} e^{if(\lambda)}\, de_{\xi , \eta}(\lambda)=\int_{\sigma(z)} e^{i\lambda}\, de_{\xi , \eta}(\lambda)=<e^{iz}\xi , \eta >\]
Therefore, we have $e^{if(z)}=e^{iz}$. Moreover, note that $f(z) \in \M _{h}$ with $\| f(z)\|\leq \pi$. Now we assert that the curve $\delta_1(t)=e^{itf(z)}$ is shorter than the curve $\delta$. Recall that for a $\tau$-measurable operator $z$ the $t$-th generalized s-number $\mu_t(z)$ is defined by
\[ \mu_t(z)=\inf \{ \, \| ze\| \, : \, e \text{ is a projection in } \M, \, \tau(1-e)\leq t \, \} \]
We shall use the following facts (see \cite{kosakifack}):
\begin{itemize}
	\item Since the map $t \mapsto \mu_t(z)$ is non-increasing,  continuous from the right and satisfies $\displaystyle{ \lim_{t\downarrow 0} \mu_t(z)}=\|z\|$, there exists $\epsilon >0$ such that $\mu_t(z)>\pi$ for all $t \in (0,\epsilon)$.
	\item It is apparent that $\mu_t(f(z))\leq \pi$, $t>0$, since $\|f(z)\|\leq \pi$.
	\item Note that $|f(t)| \leq |t|$, for all $t \in \mathbb{R}$. Then, we have $|f(z)| \leq |z|$, which implies $\mu_t(f(z)) \leq \mu_t(z)$.
\end{itemize}
Therefore,
\begin{align*}
L_p(\delta_1)^p & = \| f(z) \|_p ^p = \int_0 ^1 \mu_t(f(z)) ^p \, dt \leq \int_0 ^{\epsilon} \mu_t(f(z)) ^p \, dt + \int_{\epsilon} ^1 \mu_t(z) ^p \, dt \\
& <  \int_0 ^{\epsilon} \mu_t(z) ^p \, dt + \int_{\epsilon} ^1 \mu_t(z) ^p \, dt = \| z \|_p ^p= L_p(\delta)^p,
\end{align*}
and our claim follows.
\end{rem}

\begin{rem}\label{simbolo}
The previous remark in fact shows that, for a one-parameter group $e^{tz}$ to be minimizing, the symbol $z$ has to be bounded (and $\|z\|\le \pi$). This is due to the fact that the Borel functional calculus can be computed also for elements $z\in L^p(\M)_h$, and one obtains a shorter path by trimming the unessential parts of $z$, obtaining a shorter curve joining $1$ and $e^{iz}$.
\end{rem}

\section{Rectifiable distance in $\o$}\label{section4}

Let $x\in \o$, let $G,\g$ indicate the isotropy group and algebra of $x$ respectively. Let $\gamma\subset \o$ such that $\gamma(0)=x$. Let $\Gamma,\Lambda$ be smooth lifts of $\gamma\in \o$. Then $\alpha=\Gamma^*\Lambda\in G$, thus $\dot{\alpha}\in \alpha \g$, namely $
\dot{\Gamma^*}\Lambda+\Gamma^* \dot{\Lambda}=\Gamma^*\Lambda z$ for some $z\in \g$. Then, since $\dot{\Gamma^*}=-\Gamma^*\dot{\Gamma}\Gamma^*$, multiplying by $\Lambda^*\Gamma$ to the the right yields
\begin{equation}\label{invariante}
-\Gamma^*\dot{\Gamma}+\Gamma^*\dot{\Lambda}\Lambda^*\Gamma=Ad_{\Gamma^*\Lambda}z=\tilde{z}\in \g.
\end{equation}
Since $\|\Gamma^*\dot{\Gamma}-Q(\Gamma^*\dot{\Gamma})\|_p\le \|\Gamma^*\dot{\Gamma}+s\|_p$ for any $s\in \g$, if  we put $s=-\Gamma^*\dot{\Gamma}+\Gamma^*\dot{\Lambda}\Lambda^*\Gamma-Ad_{\Gamma^*\Lambda} Q(\Lambda^*\dot{\Lambda} )$ we obtain
$$
\|\Gamma^*\dot{\Gamma}-Q(\Gamma^*\dot{\Gamma})\|_p\le \|\dot{\Lambda}\Lambda^*-\Lambda Q(\Lambda^*\dot{\Lambda})\Lambda^*\|_p=\|\Lambda^*\dot{\Lambda}-Q(\Lambda^*\dot{\Lambda})\|_p.
$$
Since the reversed inequality also holds in the above remark, we have a natural definition of quotient metric for smooth curves $\gamma\in \o$:
\begin{defi}
Let $\gamma\in \o$ be a piecewise smooth curve. The $p$-length of $\gamma$ is defined as follows:
$$
L_{\o,p}(\gamma)=\int_0^1 \|\dot{\gamma}\|_{\gamma,p}\,dt, \mbox{ where }\|\dot{\gamma}\|_{\gamma,p}=\|\Gamma^*\dot{\Gamma}-Q(\Gamma^*\dot{\Gamma})\|_p
$$
for any smooth lift $\Gamma\in \um$ such that $\Gamma(0)=1$. The rectifiable distance in $\o$ is defined accordingly,
$$
d_{\o,p}(u\cdot x,v\cdot x)=\inf\{ L_{\o,p}(\gamma): \,\gamma\subset\o,\; \gamma(0)=u\cdot x,\,\gamma(1)=v\cdot x\},
$$
where the curves $\gamma$ considered are piecewise smooth.
\end{defi}

\begin{rem}\label{obs dist}
All but one of the properties of a distance function are satisfied by $d_{\o,p}$ trivially. The point is to establish if $d_{\o,p}(u\cdot x,v\cdot x)=0$ implies that $u\cdot x=v\cdot x$. This can not be solved as in Riemannian or Finsler geometry where the existence of normal neighborhoods is guaranteed. See Corollary \ref{dacero} below for the proof.
\end{rem}

\begin{rem}
The previous definition can be adapted for any norm quotient norm. For instance, for the uniform norm, let $\gamma:[0,1]\to\um$ be piecewise smooth, and $\Gamma:[0,1]\to \um$ is a piecewise smooth lift of $\gamma$ with $\Gamma(0)=1$, then put
$$
\|\dot{\gamma}(t)\|_{\gamma(t),\infty}=\inf_{z\in \g_x}\|\Gamma^*(t)\dot{\Gamma}(t)+z\|.
$$
The computation in (\ref{invariante}) shows that if $\Lambda$ is any other smooth lift of $\gamma$, then
$$
\inf_{z\in \g_x}\|\Gamma^*(t)\dot{\Gamma}(t)+z\|\le \|\Lambda^*(t)\dot{\Lambda}(t)+s\|,
$$
for any $s\in \g_x$. Thus, the quotient speed is well defined in this case also. The rectifiable length $L_{\o,\infty}$ and distance $d_{\o,\infty}$ are defined accordingly.
\end{rem}

\subsection{Almost isometric liftings}

We begin this section with an elementary observation, which will be used to obtain liftings of curves in $\o$. We assume that $p>1$.

\begin{lem}\label{aproxi}
Let $x\in {\cal O}$  and  $Q=Q_{{\cal G}_x}$ be the best approximant projection in $L^p(\M)_{sh}$. Let $\Gamma\subset\um$ be a piecewise smooth curve parametrized in the  interval $[0,1]$, and let $\epsilon>0$. Then there exists a polygonal curve $w_{\epsilon}:[0,1]\to {\cal G}_x$ such that $\|w_{\epsilon}(t)+Q(\Gamma^*\dot{\Gamma})(t)\|_p<\epsilon$ for any $t\in [0,1]$.
\end{lem}
\begin{proof}
Let $\alpha(t)=-Q(\Gamma^*(t)\dot{\Gamma}(t))$. Then, since $\Gamma$ is smooth for the uniform topology in $\M$, both $\Gamma$ and $\dot{\Gamma}$ are continuous for the $p$-norm, hence $\alpha:[0,1]\to L_p(\M)$ is continuous (since $Q$ is continuous). The curve $\alpha$ has its image contained in $\overline{\g_x}^p$. Then one can find a polygonal curve $w_{\epsilon}\subset \g_x$ as claimed since $\g_x$ is dense in $\overline{\g_x}^p$, as follows: split the interval $[0,1]$ in $n$ pieces $\{I_k\}_{k=1\dots n}$ in order to obtain
$$
\|\alpha(t)-\alpha(s)\|_p<\epsilon/3=\delta
$$
if $|s-t|\in I_k$, the partition given by $0=t_1<t_1<\cdots<t_n=1$, and put $\alpha_k=\alpha(t_k)$. Let $\{w_k\}_{k=1\dots n}\subset \g_x$ such that $\|\alpha_k-w_k\|_p<\delta$, and let $w_{\epsilon}(t)$ stand for the polygonal in $\g_x$ joining the points $w_k$ in their given order. Now, if $t\in I_k$, then
\begin{eqnarray}
\|w_k-w_{\epsilon}(t)\|_p& \le & \frac12 \|w_k-w_{k+1}\|_p \nonumber\\
&\le& \frac12\left[ \|w_k-\alpha_k\|_p+\|\alpha_k-\alpha_{k+1}\|_p+\|\alpha_{k+1}-w_{k+1}\|\right]\nonumber\\
&<&\frac{1}{2} [\delta+\delta+\delta]=\frac32 \delta,\nonumber
\end{eqnarray}
and hence
\begin{eqnarray}
\|\alpha(t)-w_{\epsilon}(t)\|_p&\le& \|\alpha(t)-\alpha_k\|_p+\|\alpha_k-w_k\|_p+\|w_k-w(t)\|_p\nonumber\\
&<& \frac12 \delta+ \delta+\frac32\delta=3\delta=\epsilon.\nonumber
\end{eqnarray}
\end{proof}


We collect in the following theorem some facts; some of them are elementary while others are adaptations and improvements of results in \cite{alr}.

\begin{teo}\label{levantada}
\begin{enumerate}
\item Let $k\ge 1$, $w\in \M$ with $\|w\|<\frac{\pi}{2}$. Then $T=1+\frac{(\ad w)^2}{4k^2\pi^2}$ is invertible in ${\cal B}(\M)$ and $\|T^{-1}\|\le (1-\frac{\|w\|^2}{k^2\pi^2} )^{-1}$.
\item Consider  $\displaystyle g(r)=r\sin(r)^{-1}$  with $g(0)=1$. Then $g:[0,\pi)\to \mathbb R$ is positive and increasing, and from the Weierstrass expansion of $\sin(z)$ we obtain
$$
g(z)=\prod_{k\ge 1}(1-\frac{z^2}{k^2\pi^2} )^{-1},
$$
for any $z$ such that $|z|<\pi$. Let $F(z)=(e^z-1)z^{-1}$, $w\in \M$ with $\|w\|<\frac{\pi}{2}$. Then $\|F( \ad w)^{-1}\|\le \,g(\|w\|).$
\item Let $x\in {\cal O}$, $w\subset{\cal G}_x$ a piecewise smooth curve parametrized in the  interval $[0,1]$. If $G(z)=z^{-1}(1-e^{-z})$, then there exists a piecewise smooth curve $z:[0,1]\to {\cal G}_x$ with $z(0)=0$ such that $G(\ad z)\dot{z}=w.$ If $u=e^z$, then $u:[0,1]\to G_x\subset\um$ obeys the differential equation $\dot{u}u^*=w.$
\end{enumerate}
\end{teo}
\begin{proof}
\begin{enumerate}
\item Since $\|\ad w\|\le 2\|w\|<\pi$, the map $T$ is invertible and its inverse can be computed with the Neumann series.
\item The Weierstrass expansion of $F$ is given by $F(z)=\prod_{k\ge 1}(1+\frac{z^2}{4k^2\pi^2})$ where the product converges uniformly on compact sets to $F$. Then $F(\ad w)$ is invertible since $\|ad(w)\|<\pi$ and
$$
F(\ad w)^{-1}= \prod_{k\ge 1}
\left(1+\frac{(\ad(w))^2}{4k^2\pi^2}\right)^{-1}.
$$
Hence
$$
\|F(\ad w)^{-1}\|\le \prod_{k\ge 1} \left(1-\frac{\|w\|^2}{k^2\pi^2}\right)^{-1}=g(\|w\|)
$$
by the previous item.
\item Assume first that $w$ is smooth in the whole $[0,1]$. Let $R_0=\max\limits_{t\in \overline{J}}\|w(t)\|$, where $J$ is an open set containing $[0,1]$ where $w$ is differentiable. Let $0<R<\frac{\pi}{2}$. Then if $x\in {\cal G
}_x\cap B(0,R)$, the operator $G(\ad x)$ is invertible, and its inverse is analytic and can be written as a power series in $\ad x$, hence $G(\ad x)^{-1}:{\cal G}_x\to{\cal G}_x$ because ${\cal G}_x$ is a Banach-Lie algebra. Moreover, since $g$ is increasing,
$$
\|G(\ad x)^{-1}\|\le g(\|x\|)\le g(R).
$$
Let $f:J\times B(0,R)\cap{\cal G}_x\to{\cal G}_x$ be given by $f(t,x)=G(\ad x)^{-1}w(t).$ Then $f$ is continuous since $w$ and $G^{-1}$ are continuous, moreover
$$
\|f(t,x)\|\le \|G(\ad x)^{-1}\| \, \|w\|\le g(R)R_0=L
$$
by the previous item. Now since $H(\ad x)=G(\ad x)^{-1}$ is analytic in the ball $\|x\|< \frac{\pi}{2}$, we have
$$
\|H(\ad x)-H(\ad y)\|\le C(R) \|\ad x-\ad y\|\le 2C(R)\|x-y\|
$$
where $C(R)$ is a bound for $H'$ in $|z|\le R$. Then
$$
\|f(t,x)-f(t,y)\|\le 4C(R)R_0\|x-y\|=K\|x-y\|.
$$
Then $f$ satisfies a Lipschitz condition, uniformly respect to $t\in J$, hence by the standard theorem of existence for ordinary differential equations \cite[Proposition 1.1 Ch. IV]{lang}, there exists a continuous solution $z_0: (-b,b)\times B(0,R/4)\to {\cal G}_x\cap B(0,R)$ of the integral equation
$$
z(t)=\int_0^t f(s,z(s))\,ds
$$
with $z_0(0)=0$. Here $b$ is any real number such that $0<b<\frac{R}{4LK}=\frac{\sin(R)}{32C(R)R_0^2}$. Note that $z_0$ is in fact smooth. Differentiating both sides and multiplying by $F(\ad z(t))$ gives the equation stated. We have proved so far that the equation
$$
G(\ad z)\dot{z}=w
$$
has a local solution defined around zero. By a standard argument, it follows that one can find a piecewise smooth solution defined on the whole interval $[0,1]$: let $N\in \mathbb N$ such that $\frac{1}{N}<b$ and let $t_k=\frac{k}{N}$. Then $[t_k,t_{k+1}]$ ($k=0,1,\cdots N$) is a partition of $[0,1]$ such that the integral equation
$$
z(t)=\int_{t_k}^{t_{k+1}} f(s,z(s))\,ds
$$
with the initial conditions $z_0(0)=0$, $z_k(t_k)=z_{k-1}(t_k)$ for $k\ge 1$, has a solution $z_k:[t_k,t_{k+1}]\to \cal G$. Then the curve $z_1 \sharp z_2\sharp\cdots\sharp z_N$ is a piecewise smooth solution of the equation stated in the whole $[0,1]$. If $w$ is piecewise smooth instead of smooth, one might replace the argument above for a similar argument in each of the intervals where $w$ is smooth, and use the continuity of $w$ to state the boundary conditions for $z$. If $u(t)=e^{z(t)}$, then
$$
\dot{u}(t)=exp_{*z(t)}(\dot{z}(t))=e^{z(t)}F(\ad z(t))\dot{z}(t)
$$
by Lemma \ref{expo}. Then
$$
\dot{u}u^*=e^z F(\ad z)\dot{z} e^{-z}=G(\ad z)\dot{z}=w.
$$
\end{enumerate}
\end{proof}

Note that the general theory ensures the existence of piecewise smooth liftings in $\um$ of smooth curves in $\o$, due to the fact that for any fixed $x\in\o$, the map
$$
\pi_x:\um\to \o ,\qquad \pi_x(u)=u\cdot x,
$$
is a submersion.

\begin{teo}\label{lift}
Let $\gamma\subset {\cal O}$ be a smooth curve defined in an interval containing $[0,1]$ such that $\gamma(0)=x$. Then, for any $\epsilon>0$, $\gamma$ admits a smooth lift $\beta_{\epsilon}\subset\um$ (that is $\beta_{\epsilon}\cdot x=\gamma$) such that $L_p(\beta_{\epsilon})< L_{\o,p}(\gamma)+\epsilon$. We shall call such $\beta_{\epsilon}$ an $\epsilon$-\textbf{isometric lift} of $\gamma$.
\end{teo}
\begin{proof}
Let $\Gamma\in\um$ be any piecewise smooth lift of $\gamma$, defined in an interval containing $[0,1]$, and let $w_{\epsilon}:[0,1]\to \g_x$ be as in Lemma \ref{aproxi}. Note that $w_{\epsilon}$, being a polygonal, is continuous for the uniform topology of $\M$. By item 3 of Theorem \ref{levantada}, there exists a piecewise smooth curve $u:[0,1]\to G_x$ with $u(0)=1$ such that $\dot{u}u^*=w_{\epsilon}$. Now consider $\beta_{\epsilon}=\Gamma u$. Then $\beta_{\epsilon}$ is clearly a lift of $\gamma$ with
$$
\dot{\beta_{\epsilon}}=\dot{\Gamma}u+\Gamma\dot{u}=\Gamma(\Gamma^*\dot{\Gamma}+w_{\epsilon})u.
$$
Hence $L_p(\beta_{\epsilon})<L_{\o,p}(\gamma)+\epsilon$ because
$$
\|\dot\beta_{\epsilon}\|_p = \|\Gamma^*\dot{\Gamma}+w_{\epsilon}\|_p\le \|\Gamma^*\dot{\Gamma}-Q(\Gamma^*\dot{\Gamma})\|_p+\|Q(\Gamma^*\dot{\Gamma})+w_{\epsilon}\|_p< \|\dot{\gamma}\|_{\gamma}+\epsilon.
$$

\end{proof}

With the last theorem at hand, we can prove the fundamental result that the rectifiable distance in $\o$ can be computed as the infima of lengths of rectifiable curves in $\um$ joining the corresponding fibers.

\begin{coro}\label{aprox}
Let $u,v \in \um$, $x\in \o$. Then
$$
d_{\o,p}(u\cdot x,v\cdot x)= \inf\{ \, L_p(\Gamma) \, : \, \Gamma \subset \um, \, \, \Gamma (0)\cdot x=u\cdot x \, \, \mbox{ and }\, \, \Gamma (1)\cdot x=v \cdot x \, \},
$$
where the curves $\Gamma$ considered are piecewise smooth.
\end{coro}
\begin{proof}
It suffices to check the assertion for $u=1$.  Let $d=d_{\o,p}(x, v\cdot x)$ and $D=\inf\{ \, L_p(\Gamma) \, : \, \Gamma \subset \um, \, \, \Gamma (0)\cdot x=u\cdot x \,  \mbox{ and }\, \, \Gamma (1)\cdot x=v\cdot x \, \}$. Let $\Gamma\subset\um$ such that $\Gamma(0)\cdot x=x$, $\Gamma(1)\cdot x=v\cdot x$. Then, if $\gamma=\Gamma\cdot x\subset\o$, we have $\gamma(0)=x$ and $\gamma(1)=v\cdot x$ and also
$$
\|\dot{\gamma}\|_{\gamma}=\|\Gamma^*\dot{\Gamma}-Q(\Gamma^*\dot{\Gamma})\|_p\le \|\Gamma^*\dot{\Gamma}\|_p=\|\dot{\Gamma}\|_p,
$$
thus $d\le L_{\o,p}(\gamma)\le L_p(\Gamma)$. It follows that $d\le D$. On the other hand, let $\gamma\subset \um$ joining $x$ to $v\cdot x$ such that $L_{\o,p}(\gamma)<d+\epsilon$. By the previous theorem, there exists an $\epsilon$-isometric lift $\beta_{\epsilon}$ of $\gamma$; note that $\beta_{\epsilon}(0)\cdot x=1\cdot x=x$ and $\beta_{\epsilon}(1)\cdot x=\gamma(1)=v\cdot x$. Thus
$$
D\le L_p(\beta_{\epsilon})< L_{\o,p}(\gamma)+\epsilon<d+2\epsilon.
$$
\end{proof}

\begin{coro}\label{dacero}
Let $x \in \o$. Then the quantity $d_{\o,p}$ defines a distance in $\o$ whenever $G_x$ is a closed subgroup of $\um$ in the $p$-norm.
\end{coro}
\begin{proof}
As mentioned earlier, it remained to check that $d_{\o,p}(x,y)=0$ implies $x=y$. Let $y=v\cdot x$ for some $v\in \um$, assume that $d_{\o,p}(x,y)=0$. Then by the previous corollary for any $\epsilon>0$ there exists a curve $\Gamma\in \um$ such that $\Gamma(0)=1$, $\Gamma(1)\cdot x=v\cdot x$ and the assumption implies $L_p(\Gamma)<\epsilon$. Since $\Gamma(1)\in vG_x$, then $d_p(1,vG_x)\le L_p(\Gamma)<\epsilon$. Since $\epsilon$ is arbitrary, then $1\in vG_x$, or equivalently $v^*\in G_x$, and then $v\in G_x$, showing that $y=x$.
\end{proof}

\begin{rem}
We point out that when $G_x$ is the unitary group of a von Neumann subalgebra of $\M$, then $G_x$ is $p$-norm closed in $\um$.
\end{rem}

\subsection{Two metrics in the space $\o$}

For homogeneous spaces $G/H$ of metrizable topological groups $G$, there is one distinguished metric that can be introduced. It is the quotient metric induced by the distance among classes $g H$ in the original group. We recall the following result (see for instance \cite[p.109]{takesaki3}):

\begin{lem}\label{take}
Let $G$ be a metrizable topological group, and $H$ be a closed subgroup. If $d$ is a complete distance function on $G$ inducing the topology of $G$, and if $d$ is invariant under right translation by $H$, i.e. $d(xh,yh)=d(x,y)$ for any $x,y \in G$ and $h \in H$, then the left coset space $G/H=\{ \, xH \, : \, x \in G \, \}$ is a complete metric space under the metric $\dot{d}$ given by
\[ \dot{d}(xH,yH)= \inf \{ \, d(xh,yk)\, : \, h,k \in H \, \}.  \]
\end{lem}

\medskip

Moreover, the distance $\dot{d}$ metrizes the quotient topology of groups. Let us observe how Lemma \ref{take} applies to our situation. We take $G=\um$, and for fixed $x\in \o$, we take $H=G_x$.

\begin{teo}\label{iguales}
Let $x\in \o$, $u,v\in \um$, and let
$$
\dot{d}_p(u\cdot x,v\cdot x)=\inf \{ \, d_p(uw_1,vw_2)\, : \, w_i  \in G_x \, \}.
$$
Then if $p>1$, $\dot{d}_p=d_{\o,p}$. In particular, if  $G_x$ is a closed subgroup of $\um$ in the $p$-norm, then $(\o,d_{\o,p})$ is complete, and the induced topology matches the quotient topology of $\o\simeq (\um,d_p)/G_x$.
\end{teo}
\begin{proof}
First we show that $\dot{d}_p \leq d_{\o,p}$. By Corollary \ref{aprox}, for each $\epsilon >0$, there exists a curve $\Gamma \subseteq \um$ satisfying $\Gamma(0)=uw_0$, $\Gamma(1)=vw_1$, $w_i \in G_x$ and $L_p(\Gamma) < d_{\o,p}(u\cdot x,v\cdot x) + \epsilon$. Therefore,
$$
\dot{d}_p(u\cdot x,v\cdot x) \leq d_p(uw_0,vw_1) \leq L_p(\Gamma) < d_{\o,p}(u\cdot x,v\cdot x) + \epsilon.
$$
Since $\epsilon$ is arbitrary, our claim follows. Conversely, given $\epsilon >0$, there exist $w_i \in G_x$, $i=1,2$ such that $d_p(uw_1, vw_2) < \dot{d}_p(u\cdot x,v\cdot x) + \epsilon$, and there exists $\Gamma\subset\um$ such that $\Gamma(0)=uw_1$, $\Gamma(1)=uw_2$, and $L_p(\Gamma)<d_p(uw_1,uw_2)+\epsilon$. Then
$$
d_{\o,p}(u\cdot x,v\cdot x) \leq L_{\o,p}(\Gamma \cdot x) \leq L_p (\Gamma) < d_p(uw_1, v w_2)+\epsilon< \dot{d}_p(u\cdot x,v\cdot x) + 2\epsilon,
$$
showing that the reversed inequality also holds. Finally, since $(\um,d_p)$ is complete, the last claim follows from the previous lemma.
\end{proof}

\subsection{The path metric space $\o$}

The space  $\o$ is a \textit{path metric space} for any $p\ge 1$, or in the terminology introduced in \cite{gromov} by M. Gromov, $\o$ is a \textit{space de longeur}. That is, the distance $d_{\o,p}$ ($=\dot{d}_p$ if $p>1$) among pairs of points in the space $\o$ matches the infimum of the length of the rectifiable paths joining the points. The rectifiable lenght of paths $\gamma:[0,1]\to \o$ is defined as
$$
\ell_p(\gamma)=\sup\limits_{\{t_i\}} \sum\limits_{i=0}^{n-1} d_{\o,p}(\gamma(t_i),\gamma(t_{i+1})),
$$
where $\{t_i\}$ is any finite partition of the interval $[0,1]\in\mathbb R$. The rectifiable distance $d_{\ell,p}$ among $x,u\cdot x\in \o$ is the infimum of the length of the rectifiable paths $\gamma$ joining them, when the length is measured as above. It is straightforward to see that
$$
d_{\ell,p}\ge d_{\o,p}.
$$
Indeed, for given $\gamma$ joining fixed endpoints $x,y\in \o$, consider the trivial partition $t_0=0,t_1=1$. Thus $\ell_p(\gamma)\ge d_{\o,p}(x,y)$, and taking the infimum over rectifiable paths $\gamma$ gives the result. It is well-known, at least in the finite dimensional setting, that both metrics do agree. Since the proof is elementary and we could not find a suitable reference, we include it.

\begin{prop}
If $\gamma$ is a piecewise $C^1$ curve in $\o$, then $L_{\o,p}(\gamma)\ge \ell_p(\gamma)$. If $x, u\cdot x\in \o$, then $d_{\ell,p}(x,u\cdot x)=d_{\o,p}(x, u\cdot x)$.
\end{prop}
\begin{proof}
Let $\{t_i\}_{i=0\cdots n-1}$ be a partition of $[0,1]$. Then
$$
L_{\o,p}(\gamma)=\sum\limits_{i=0}^{n-1} L_{\o,p}(\gamma|_{[t_i,t_{i+1}]})\ge\sum\limits_{i=0}^{n-1} d_{\o,p}(\gamma(t_i),\gamma(t_{i+1})).
$$
If we pick a partition such that $\ell_p(\gamma)<\sum\limits_{i=0}^{n-1} d_{\o,p}(\gamma(t_i),\gamma(t_{i+1}))+\epsilon$, the first claim follows. Now let $\epsilon>0$, and let $\Gamma:[0,1]\to\um$ be piecewise smooth such joining $x,u\cdot x$ such that $L_p(\Gamma)\le d_{\o,p}(x,u\cdot x)+\epsilon$. Then by the previous assertion
$$
d_{\ell,p}(x,u\cdot x)\le \ell_p(\Gamma \cdot x)\le L_{\o,p}(\Gamma\cdot x) \le L_p(\Gamma)\le d_{\o,p}(x,u\cdot x)+\epsilon.
$$
Thus $d_{\ell,p}(x,u\cdot x)\le  d_{\o,p}(x,u\cdot x)$ and since the other inequality always holds, we have the second claim.
\end{proof}

\section{Minimality of geodesics in $\o$}\label{section5}

Recall that the induced norm in the tangent spaces is not complete. Therefore  in the case of $p=2$ the classical theory of Riemann-Hilbert manifolds is not available, so it makes sense to ask about the local minimality of the geodesics of the Levi-Civita connection. In \cite{chiumiento} was given an abstract sufficient condition in order that these geodesics are locally minimal. In this section we shall prove a partial result toward the minimality of geodesics for $p$ even, under the hypothesis specified below.

Our argument on minimality will consist in comparing the lengths of the liftings of curves in $\o$ to the unitary group $\um$. For the case $p=2$ this technique is based on the following result. Let $s_x:(T\o)_x\to  \f_{x}\subset \beah$ stand for the isometric orthogonal projection to $\f_{x}$, the orthogonal supplement of $\g_x$ in $\beah$. Let $\gamma(t)$, $t\in[0,1]$ be a smooth curve in $\o$, with $\gamma(0)=x$,  and let $\Gamma$ be its horizontal lifting, i.e. the unique solution of the differential equation
$$
\left\{
\begin{array}{l}
\dot{\Gamma}= s_\gamma(\dot{\gamma}) \Gamma \\
\Gamma(0)=  1 .
\end{array}
\right.
$$
Then $L_2(\Gamma)=L_2(\gamma)$. Indeed, $\|\dot{\gamma}\|_\gamma=\|s_\gamma(\dot{\gamma})\|_2=\|\Gamma^* \dot{\Gamma}\|_2=\|\dot{\Gamma}\|_2$, and the claim follows. This result shall not be needed, we include it here to mark the breach between the  case $p=2$ and the case $p>2$.   Let us state the following definition.

\begin{defi}
Let $\o=\um\cdot x$ ($x\in \o)$ be an homogeneous space, let $G_x\subset \um$ stand for the isotropy subgroup. We say that $G_x$ is locally exponential in $\um$, if  there exist $\epsilon_{\o},\delta_{\o}>0$ such that $\|u-1\|<\epsilon_{\o}$ and $u\in G_x$ implies that there exists $z\in \g_x$ with $\|z\|<\delta_{\o}$ and $e^z=u$. This is equivalent to ask for $G_x$ to be a (topological) submanifold of $\um$ in the uniform norm.

If $u\in G_x$ implies that there exists $z\in  \g_x$ with $e^z=u$, we say that $G_x$ is an exponential subgroup of $\um$.
\end{defi}

\noi Throughout this section we assume that the isotropy group $G_x$ is an exponential subgroup. Apparently, if this holds for a given $x\in\o$, then it holds for any $u \cdot x \in \o$ (since the groups $G_{u \cdot x}$ and $G_x$ are conjugate by an inner automorphism). This property implies in particular, that $G_x$ is geodesically convex: given any pair of elements $v_1,v_2 \in G_x$ with $\|v_1-v_2\|<\epsilon_{\o}$, then there exists a geodesic of $\um$, which lies inside $G_x$, and joins $v_1$ and $v_2$. The results in this section can be extended to locally exponential groups in the obvious fashion, but we prefer to state them assuming that $G_x$ is exponential since the arguments become more clear this way.

Assume that $x$ and $y$ are connected by the geodesic $\gamma(t)=e^{tz}\cdot x$ in $\o$, with  $z\in\M_{sh}$ a minimal lifting, i.e. $Q(z)=0$. It is unclear if this curve is short for the $p$-metric in $\o$. We state next a partial result in that direction.

One requirement of the proof is that $Q$ maps bounded operators into bounded operators, and moreover, that it is uniformly bounded. There are a collection of examples in the following section with these properties. We have shown in Theorem \ref{lift} that for $p>2$ we can obtain almost isometric lifts of curves in ${\cal O}$. In case $Q$ is uniformly bounded, one can sharpen this result to obtain isometric lifts.

\noi Now we state our  result on minimality of curves. It is assumed that $G_x$ is an exponential subgroup of $\um$, and that $Q$ maps bounded elements of $\M$ into bounded elements of $\M$.

\begin{teo}\label{geod}
Let $p$ be a positive even number, $x\in\o$, and assume that there exists a constant $K_{\o,p}$ such that $\|Q(x)\|\le K_{\o,p}\|x\|$ for any $x\in \M_{sh}$. If $z\in \M_{sh}$, $\|z\|<\frac{\pi}{3}$ and $Q(z)=0$, then the curve
$$
\delta(t)=e^{t z}\cdot x,
$$
which verifies $\delta(0)=x$ and $\dot{\delta}(0)=X=\pi_*(z)\in T_x\o$, and has length $L_{\o,p}=\|z\|_p$, is shorter for the $p$-metric than any other smooth curve $\gamma\subset\o$ joining $x$ to $e^z\cdot x$, provided $L_{\o,\infty}(\gamma)<\varepsilon$, where
$$
\varepsilon=\varepsilon(\o,p)=\frac{\sqrt{2}-1}{C_{\o}(1+K_{\o,p})}
$$
and $C_{\o}$ is a constant given by the smooth structure as in Remark \ref{parte}.

Moreover, the curve $\delta$ is unique in the sense that if $\gamma\subset {\cal O}$ is another curve joining $x$ to $e^z\cdot x$ of length $\|z\|_p$, such that $L_{\o,\infty}(\gamma)<\varepsilon$, then $\gamma(t)=e^{tz}\cdot x$.
\end{teo}
\begin{proof}
Let $\gamma$ be a smooth curve in $\o$ with $\gamma(0)=x$ and $\gamma(1)=e^{z}\cdot x$, and assume that $L_{\o,\infty}(\gamma)<\varepsilon$. Since $\|z\|<\frac{\pi}{3}$, we have $
\|e^z-1\|<1<\sqrt{2}$. On the other hand, if $\Gamma$ is a smooth lift of $\gamma$ with $\Gamma(0)=1$, by the assumption on $Q$ we can consider the differential equation in $\M$ given by
$$
G(\ad x)\dot{x}=-Q(\Gamma^*\dot{\Gamma})
$$
as in Theorem \ref{levantada}. It has a unique solution $x(t)\in\g_x$ such that $x(0)=0$, and if $u=e^x$, then $u:[0,1]\to G_x\subset\um$ obeys the differential equation $\dot{u}u^*=-Q(\Gamma^*\dot{\Gamma})$. Thus in this case, $\|\dot{u}\|\le K_{\o}\|\dot{\Gamma}\|.$ Hence, if $\beta=\Gamma u$, then $\beta$ is an isometric lift of $\gamma$ and moreover
$$
\|\beta(1)-1\|\le \|\Gamma(1)-1\|+\|u(1)-1\|\le \int_0^1 \|\Gamma^*\dot{\Gamma}\|+\int_0^1 \|\dot{u}\|\le (1+K_{\o}) \int_0^1 \|\Gamma^*\dot{\Gamma}\|.
$$
Thus, if $C_{\o}$ is a constant such that $\|s_y(V)\|\le C_{\o}\|V\|_y$ for any $y\in \o$ as in Remark \ref{parte}, then
$$
\|\beta(1)-1\|\le C_{\o} (1+K_{\o})\int_0^1 \|\dot{\gamma}\|_{\gamma}=C_{\o} (1+K_{\o})L_{\o,\infty}(\gamma)<\sqrt{2}-1<\sqrt{2}-\|e^z-1\|.
$$
Let $w\in\M_{sh}$ such that $\|w\|\le \pi$ and $e^w=\beta(1)$. Then Theorem \ref{minimalidadunitarios} applies. Note that the curve $\mu(t)=e^{tz}$ is an isometric lift for $\delta$. Let $\nu(t)=e^{z}e^{ty}$ be the minimal geodesic of $\um$, lying inside $e^{z} G_x$ (i.e. $y\in\g_x$), connecting $e^z$ to $e^w$, which exists due to the fact that $G_x$ is an exponential subgroup. Then by item 5 of Theorem \ref{minimalidadunitarios}, the map $f(s)=d_p^p(1,\nu(s))$ is convex. Note that $f'(0)=(-1)^{p/2} Tr(z^{p-1}y)$, which vanishes by Lemma \ref{perpep}, because $z$ is a minimal lift. Then
$$
L_p(\mu)^p=d_p^p(1,\nu(0))=f(0)\le f(1)=d_p^p(1,\nu(1))=\|w\|_p^p\le L_p(\beta)^p.
$$
Hence
$$
\|z\|_p=L_{\o,p}(\delta)\le L_p(\beta)=L_{\o,p}(\gamma).
$$

\smallskip

\noi If $L_{\o,p}(\gamma)=\|z\|_p$ (i.e. if $\gamma$ is also short), then
$$
f(0)=\|z\|_p\le f(1)=\|w\|_p\le L_p(\beta)=L_{\o,p}(\gamma)=\|z\|_p=f(0),
$$
so $f(1)=f(0)$, which forces $z=w$ because $f$ is strictly convex. In particular $\beta(1)=e^z$ and $L_p(\beta)=L_p(\mu)=\|z\|_p$. Since $\|z\|<\pi/2<\pi$, the curve $\mu$ is the unique short geodesic joining $1$ to $e^{z}$ in $\um$, and then $\beta=\mu$, or in other words, $\gamma=\delta$.
\end{proof}

\begin{rem}\label{uniforme}
The restriction on the quotient uniform length of the curves $\gamma$ can be removed for $p=2$  due to the existence of the supplement ${\cal F}$ given by the smooth structure of $\o$ satisfying $\M_{sh}\simeq \g\oplus {\cal F}$. The key is that the exponential map is a local diffeomorphism  between ${\cal F}$ and $\o$. In the general case $p>2$ we do not know if the exponential map is local bijection between $\g^{\perp_p}$ and $\o$ (or even if $\g^{\perp_p}\cap \M\ne \{0\}$).

However, assuming that $Q$ maps bounded elements of $\M$ into bouded elements of $\M$, it is interesting to note that the map
$$
\varphi=(1-Q_{\g})|_{{\cal F}}:{\cal F} \longrightarrow \g^{\perp_p}
$$
is a bijection. In fact, $Q(\varphi(f))=Q\circ (1-Q)(f)=0$ for any $f\in {\cal F}$ showing that $\varphi$ maps into $\g^{\perp_p}$. Secondly $\varphi(f_1)=\varphi(f_2)$ implies $f_1-f_2=Q(f_1)-Q(f_2)\in\g_x$ thus $f_1=f_2$ if $f_i\in {\cal F}$, showing that $\varphi$ is injective. In third place, if $z\in \g^{\perp_p}$, $z=z_g+z_f$ with $z_g\in\g$ and $z_f\in {\cal F}$, hence $0=Q(z)=z_g+Q(z_f)$. Thus taking $f=z_f\in {\cal F}$, one obtains $\varphi(f)=z_f-Q(z_f)=z_f+z_g=z$ showing that $\varphi$ is surjective. The inverse is given by the linear projection onto ${\cal F}$, that is
$$
\varphi^{-1}=P_{\cal F}|_{\g^{\perp_p}}:\g^{\perp_p}  \longrightarrow {\cal F}.
$$
Note that, while $\varphi$ is continuous for the $p$-norm, $\varphi^{-1}$ is continuous for the uniform norm, showing the breach between the smooth structure and the metric structure.

It is also worthwile noting that, if $Q$ is continuous for the uniform topology of $\M$, then the above maps are homeomorphisms. Moreover, it must be uniformly bounded since, if it were not,  there would exist a sequence $(x_n)_{n\ge 1}$ of elements of $\M$ such that $\|x_n\|=1$ and $\|Q(x_n)\|\ge n$. But this contradicts the fact that $Q(0)=0$, since $
\|Q(\frac{x_n}{n})\|\ge 1$, and the assumption on the continuity of $Q$ gives $Q(\frac{x_n}{n})\to 0$ in $\M$.
\end{rem}

\begin{rem}\label{uniq}
Assume that $Q$ is uniformly bounded $(\|Q(z)\|\le K_{\o,p} \|z\|$ for any $z\in \M_{sh}$). Let $B_R(0)\subset {\cal F}$ with $R$ small enough to ensure $\pi\circ \exp$ is a diffeomorphism with its image in $\o$. Consider $V_{\perp}^R=\varphi(B_R(0))$, which is open in $\g^{\perp}$ with the relative (uniform) topology, since $\varphi^{-1}=P_{\cal F}$ is continuous. Then, if $z\in V_{\perp}^R$, we have $z=\varphi(y)$ for some $y\in B_R(0)$, hence
$$
 \|z\| =\|y-Q(y)\|\leq (1+K_{\o \, , \, p})\|y\|< (1+K_{\o \, , \, p})R.
$$
Let
$$
U_{\o}^R=\pi\circ \exp(V_{\perp}^R)=\{e^{w-Q(w)}\cdot x : w\in{\cal F},\, \|w\|<R\}.
$$
Note that it is not clear whether this is an open neighborhood of $x\in \o$ or not, even if we assume that $Q$ is continuous for the uniform topology.
\end{rem}

Now we can state our theorem on minimal curves joining given endpoints in $\o$. Let $p$ be a positive even number, $x\in\o$, and assume that there exists a constant $K_{\o,p}$ such that $\|Q(x)\|\le K_{\o,p}\|x\|$ for any $x\in \M_{sh}$. Let
$$
r=\min\{R,\frac{\varepsilon}{2(1+K_{\o,p})},\frac{\pi}{3}\}
$$
where $R$ is as in the previous remark, and $\varepsilon=\varepsilon(\o,p)$ as in the previous lemma. Let $V^r_{\perp}$, $U_{\o}^r\subset \o$ be the sets defined in the previous remark.

\begin{teo}\label{geodesico}
For any $y\in U_{\o}^r$ there exists $z\in V^r_{\perp}$ such that $e^z\cdot x=y$ and
$$
\delta(t)=e^{t z}\cdot x
$$
is shorter for the $p$-metric than any other piecewise smooth curve $\gamma\subset\o$ joining $x$ to $y$, provided $\gamma\subset U_{\o}^r$.

Moreover, the curve $\delta$ is unique in the sense that if $\gamma\subset U_{\o}^r$ is another piecewise smooth curve joining $x$ to $y$ of length $\|z\|_p$ then $\gamma(t)=e^{tz}\cdot x$.
\end{teo}
\begin{proof}
The existence of such $z$ is guaranteed by Remark \ref{uniq}. Let $\gamma\subset U_{\o}^r$ be piecewise smooth, we can assume that $\gamma$ is defined in $[0,1]$. Consider a partition $\{[t_i,t_{i+1}]\}$ of $[0,1]$ in $N$ equal pieces such that $L_{\o,\infty}(\gamma|_{[t_i,t_{i+1}]})<r$. By the previous theorem, $\gamma|_{[t_0,t_1]}$ is longer than the curve
$
\delta_1(t)=e^{z_1}\cdot\gamma(t_0),
$
where $z_1\in V_{\perp}^r$ is such that $e^{z_1}\cdot\gamma(t_0)=\gamma(t_1)$. Let $\alpha \sharp\beta$ denote the path $\alpha$ followed by the path $\beta$. Then
$$
L_{\o,\infty}(\delta_1 \sharp \gamma|_{[t_1,t_2]})=L_{\o,\infty}(\delta_1)+ L_{\o,\infty}( \gamma|_{[t_1,t_2]})<\|z_1\|+r<2r.
$$
Let $z_2\in V^r_{\perp}$ such that $e^{z_2}\cdot x=\gamma(t_2)$. Then by the previous theorem, the path $\delta_2(t)=e^{tz_2}\cdot x$ is shorter than $\delta_1 \sharp \gamma|_{[t_1,t_2]}$, thus $\delta_1\sharp \gamma|_{[t_2,1]}$ is shorter that $\gamma$. Iterating this argument, one ends with a curve $\delta_N=e^{tz_N}\cdot x$, joining $x$ to $y$ in $U_{\o}^r$, which is shorter that $\gamma$. Since $\|z_N\|<r$, it must be $z_N=z$.

The uniqueness follows observing that in each step, if the length of $\gamma$ is equal to $\|z\|_p$, its restriction must have length $\|z_i\|_p$, and by the previous lemma it should match $\delta_i$.
\end{proof}

\subsection{Examples}

We give examples of homogeneous spaces where Theorems \ref{geod} and \ref{geodesico} apply. The fundamental step is to prove that the metric projection $Q$ is uniformly bounded. Up to now we do not know if it is general fact, even in the case when the Lie algebra where $Q$ projects consists of skew-hermitian operators of a von Neumann subalgebra. Therefore each example needs an ad-hoc proof of this fact. We give sketches of proofs here, full proofs can be found in \cite{eduard}.  Finally, let us observe that the examples below serve obviously as other examples for Theorem   \ref{iguales}.

\subsubsection{Finite dimensional Lie algebras}
The first  immediate example takes place when the Lie algebra $\g$  is a finite dimensional vector space. Therefore the $\| \,. \, \|_p\,$ completion of $\g$ is equal to $\g$. Hence it is trivial that the  projection $Q:\lpm_{sh} \longrightarrow \overline{\g}^p=\g$ preserves bounded elements.

\begin{lem}
Let $1<p<\infty$ and $\g$ a finite dimensional Lie algebra. Then the projection $Q$ is continuous and in particular, uniformly bounded.
\end{lem}

\begin{proof}
Since $\g$ is a finite dimensional real vector space, all the norms are equivalent. Therefore, there exist a  constant $c_p >0$ such that $c_p\,\|z\| \leq \|z\|_p \leq \|z\|$, for all $z \in \g$. Now, given $\epsilon >0$, there exists $\delta(x,\varepsilon,p)$ such that $\|x-y\|_p<\delta$ implies $\|Q(x)-Q(y)\|_p<\varepsilon$. Hence, if $\|x-y\|<\delta$, then $\|x-y\|_p<\delta$ and
$$
\|Q(x)-Q(y)\|\le c_p^{-1} \|Q(x)-Q(y)\|_p<c_p^{-1}\varepsilon.
$$
The argument at the end of Remark \ref{uniforme} establishes the uniform boundness of $Q$.
\end{proof}

\noi We describe an example where this situation arises.

\begin{ejem}
Let $v_0\in \M$ be a partial isometry of finite co-rank. Consider the set
$$\mathcal{I}_{v_0}=\{  \, v \in \M \, : \, v_0^*v_0=v^*v \, \} $$
of partial isometries in $\M$ with initial space $p$. There is a transitive action of $\um$ on $\mathcal{I}_{v_0}$ given by $u \cdot v = uv$, $u \in \um$, $v \in \mathcal{I}_{v_0}$. The set $\mathcal{I}_{v_0}$ is a $C^{\infty}$ submanifold of $\M$ in the norm topology and a homogeneous space of $\um$.
The isotropy group at $v \in \mathcal{I}_{v_0}$ of the action is
\[  \{ \, u \in \um \, : \, u v=v \, \}. \]
Therefore the Lie algebra of the above group is given by
\[ \g_v = \{ \, x \in \M_{sh} \, : \, xv =0 \, \}, \]
and the unitaries in the isotropy group can be described as
$$
u=\bigg(\begin{array}{cc}{1}&{0}\\{0}&{d}\end{array}\bigg),
$$
with $d$ a unitary operator, the group is exponential. Then Theorem \ref{geod} applies to this situation, and the curves $\delta(t)=e^{tz}v$ with minimal symbol $z$ are short among sufficiently short curves $\gamma\in \o$. This example was studied in \cite{partialiso}.
\end{ejem}

\subsubsection{Subalgebra of the center}

This example is concerned with a subalgebra $\N \subseteq \mathcal{Z} (\M)$, where $\mathcal{Z} (\M)$ is the center of $\M$. In order to show that $Q$ preserves bounded elements we have the following lemma, see \cite{eduard} for a proof and a counterexample for $p>2$, if we remove the hypothesis  that $x$ and $y$ commute.

\begin{lem}
Let $p\geq2$ an even number. Let $x,y \in \lpm$ satisfying $x \geq 0$, $y=y^*$ and $xy=yx$. Then
\[ \| x- y^+ \|_p \leq \| x - y \|_p \, , \]
where $y=y^+ - y^-$ is the Jordan decomposition.
\end{lem}

\bigskip

\noi Applying the previous lemma to a positive element $x \in \M$ and $Q(x) \in \lpn$ we obtain
\[ \| x- Q(x)^+ \|_p \leq \| x - Q(x) \|_p \,.\]
Hence by the uniqueness of $Q(x)$ it follows that $Q(x)=Q(x)^+$. In particular, it follows that  the projection $Q$ maps positive elements of $\lpm$ into positive elements of $\lpn$.

\begin{coro}
If ${\cal N}\subset {\cal Z}(\M)$, then for any $p>1$ the projection $Q$ maps bounded elements into bounded elements, moreover
$$
\|Q(z)\|\le 3 \|z\|\;\mbox{ for any }z\in \M_{sh}.
$$
\end{coro}
\begin{proof}
Let $x \in \M$ a positive element. Note that $\| x\| - Q(x)= Q(\| x \|-x)\geq 0$, then $0 \leq Q(x) \leq \| x \|$, i.e. $Q(x)$ is bounded. Let $x \in \M_{h}$, then there exists a real number $c>0$ such that $x+ c$ is positive. Since $Q(x) + c=Q(x +c)$ is bounded, it follows that $Q(x)$ is bounded, and moreover, if $x\in\M_h$ then
\begin{eqnarray}
\|Q(x)\|&=&\|Q(x+\|x\|-\|x\|)\|=\|Q(x+\|x\|)-\|x\|\|\le \|Q(x+\|x\|)\|+\|x\|\nonumber\\
&\le &\|x+\|x\|\|+\|x\|\le 3\|x\|.\nonumber
\end{eqnarray}
Replacing $x$ by $ix$ yields the result for $z\in \M_{sh}$.
\end{proof}

\begin{rem}
In case that the Lie algebra is given by antihermitic operators of a von Neumann subalgebra of $\M$ we have the bound  $C_{\o}\leq 2$. This follows because the projection $P_{{\cal F}_x}$ coincides with $I-E$, where $E$ is the unique normal conditional expectation preserving the trace onto the subalgebra.
\end{rem}

\begin{rem}
Let $U_{\o}^r\subset \o$ as in Theorem \ref{geodesico}. If $\o$ is the quotient space obtained as $\um/{\cal U}_{\cal N}$, and ${\cal N}\subset {\cal Z}(\M)$, then $U_{\o}^r$ is an open neighborhood of $x$ in $\o$. In fact,
$$
U_{\o}^r=\{e^{w-Q(w)}\cdot x : w\in{\cal F},\, \|w\|<r\}=\{e^w\cdot x : w\in{\cal F},\, \|w\|<r\},
$$
and the last set is clearly open in $\o$ by our choice of $r$.
\end{rem}

\subsubsection{Diagonal algebra in $\M \otimes M_2$}\label{alg diagonal en mat 2por2}

Let $M_2$ denote the $2 \times 2$ matrix algebra. We define a  finite trace $\hat{\tau}$ on $\M \otimes M_2$ by
\[    \hat{\tau} \big(  \, \bigg(\begin{array}{cc}{x_{11}}&{x_{12}}\\{x_{21}}&{x_{22}}\end{array}\bigg) \, \big)=\frac12\tau(x_{11} + x_{22}), \, \, \, \, \, \, \, \, \, \,  \bigg(\begin{array}{cc}{x_{11}}&{x_{12}}\\{x_{21}}&{x_{22}}\end{array}\bigg) \in  M_2 \otimes \M. \]
It is straightforward to show that $L^p(\M \otimes M_2 , \hat{\tau})= L^p(\M)\otimes M_2$.

\medskip

\noi We take the subalgebra $\N$ consisting of diagonal operator matrices, i.e.
\[  \N =   \{ \, \bigg(\begin{array}{cc}{x_{11}}&{0}\\{0}&{x_{22}}\end{array}\bigg)   \, : \, x_{11}, \, x_{22} \in \M   \, \}. \]

\noi In this example we can explicitly compute the projection $Q$. Actually, this is a consequence of the following inequality, see \cite{eduard} for a proof.

\begin{lem}\label{lema}
Let $p\geq 2$ a positive even number and $b \in \M$. Then,
\[ \bigg\| \bigg(\begin{array}{cc}{0}&{b}\\{b^*}&{0}\end{array}\bigg)  \bigg\|_p \leq \bigg\| \bigg(\begin{array}{cc}{a}&{b}\\{b^*}&{d}\end{array}\bigg) \bigg\|_p\,, \]
for all $a,b \in \M_{h}$.
\end{lem}

\noi It is plain that $Q: (L^p(\M)\otimes M_2)_h \longrightarrow \lpn_h$ is given by
\[ Q(\,\bigg(\begin{array}{cc}{x_{11}}&{x_{12}}\\{x_{21}}&{x_{22}}\end{array}\bigg)  \,)=\bigg(\begin{array}{cc}{x_{11}}&{0}\\{0}&{x_{22}}\end{array}\bigg).  \]
In particular, $Q$ preserves bounded elements and is uniformly bounded. Moreover, it is continuous in the uniform topology since it matches the unique tracial invariant conditional expectation $E$ from the algebra to the subalgebra.

\begin{ejem}
 Consider the projection in $\M \otimes M_2$ given by $e=\bigg(\begin{array}{cc}{1}&{0}\\{0}&{0}\end{array}\bigg)$.
Let $\mathcal{O}_e$ denote the unitary orbit, i.e. $ \mathcal{O}_e= \{ \, ueu^* \, : u \in \mathcal{U}_{\M \otimes M_2} \, \}.$ This example was studied in detail in \cite{otrococo}. It was proved that it is a homogeneous space of the unitary group $\mathcal{U}_{\M \otimes M_2}$ of $\M \otimes M_2$. Moreover, it was shown that the initial values problem has solution and any pair of points in this homogeneous space can be joined by a minimal curve. Despite our results are more restrictive in this particular example, we shall  show how they apply, since the techniques involved are quite different.

The isotropy group at $e$ of the natural action of $\mathcal{U}_{\M \otimes M_2}$ is given by
$$G_e= \{ u \in \mathcal{U}_{\M \otimes M_2} \, : \, ue=eu \, \}$$
The Lie algebra of this group is
\begin{align*}
\mathcal{G}_e & =\{ x \in (\M \otimes M_2)_{sh} \, : \, xe=ex \, \}
 = \{ \bigg(\begin{array}{cc}{a}&{0}\\{0}&{d}\end{array}\bigg) \, : \, a,d \in \M_{sh}\, \}.
\end{align*}
Therefore by our preceding discussion the projection $Q$ onto the Lie algebra preserves bounded elements, so our results about minimality of curves holds.
\end{ejem}

\subsubsection{Special diagonal algebra in $\M \otimes M_2$}

Consider the following subalgebra of $\M \otimes M_2$ given by
\[  \N=\{ \, \bigg(\begin{array}{cc}{x}&{0}\\{0}&{x}\end{array}\bigg) \, : \, x \in \M \, \}.\]
Let $E$ denote the unique trace-invariant (with respect to the trace $\hat{\tau}$) conditional expectation onto $\N$, i.e.
\[ E: \M \otimes M_2\longrightarrow \N, \, \, \, \, \,\, \, \, \, \,  E ( \, \bigg(\begin{array}{cc}{x_{11}}&{x_{12}}\\{x_{21}}&{x_{22}}\end{array}\bigg) \, ) =    \frac{1}{2}\bigg(\begin{array}{cc}{x_{11} + x_{22}}&{0}\\{0}&{x_{11} + x_{22}}\end{array}\bigg). \]
We denote by $E_p$ the extension of the above expectation to the corresponding non commutative $L^p$ spaces.

\begin{lem}
Let $2\leq p < \infty$, $p$ even. Then:
\[  \bigg\|  \bigg(\begin{array}{cc}{(a-c)/2}&{b}\\{b}&{(c-a)/2}\end{array}\bigg)    \bigg\|_p \leq
   \bigg\|  \bigg(\begin{array}{cc}{a}&{b}\\{b}&{c}\end{array}\bigg) + \bigg(\begin{array}{cc}{d}&{0}\\{0}&{d}\end{array}\bigg)                 \bigg\|_p   \]
for any $d \in \M$.
\end{lem}
\begin{proof}
See \cite{eduard}.
\end{proof}

If $\mathcal{L}$ stands for the following real subspace of $\M_h \otimes M_2$ given by
$$ \mathcal{L}=\{\,\bigg(\begin{array}{cc}{a}&{b}\\{b}&{c}\end{array}\bigg)\, : \, a,b,c \in \M_{h} \, \}$$
and $\mathcal{L}^p$ the respective completion with the $p$-norm. Then, it is easy to check, using the previous lemma, that $E: \mathcal{L} \longrightarrow \N$ and $E_p: \mathcal{L}^p \longrightarrow L^p(\N)$ for $p$ even, are contractive maps.

Analogous statements hold for the subspace
$$
{\cal D}=\{ \, \bigg(\begin{array}{cc}{0}&{b}\\{b^*}&{0}\end{array}\bigg) \, : \, b \in \M \, \},
$$
invoking Lemma \ref{lema}. If $p$ is even or $p=\infty$, then $Q_{\N,p}=E_p$, namely the best approximant can be obtained via the conditional expectation in ${\cal L}$. In particular $Q$ is uniformly bounded ${\cal L}$. A similar argument shows that $Q$ is uniformly bounded in ${\cal D}$.

It is not clear, and we would like to know, if $Q$ is uniformly bounded in $\M\otimes M_2$.

\bigskip

\noi 
Eduardo Chiumiento\\
e-mail: eduardo@mate.unlp.edu.ar\\
Departamento de Matem\'atica, FCE-UNLP\\
Calles 50 y 115\\
(1900) La Plata, Argentina\\
and \\
Instituto Argentino de Matem\'atica, CONICET
\\

\noi
Esteban Andruchow and Gabriel Larotonda\\
e-mails: eandruch@ungs.edu.ar, glaroton@ungs.edu.ar\\
Instituto de Ciencias\\
Universidad Nacional de General Sarmiento\\
J. M. Gutierrez 1150\\
(B1613GSX) Los Polvorines, Argentina\\
and \\
Instituto Argentino de Matem\'atica, CONICET
\\

\end{document}